\documentclass[11pt]{article}   
\usepackage{amssymb,amscd,latexsym}   
\usepackage{amsmath}
\usepackage{amsthm}
\usepackage{amsfonts}
\usepackage{mathdots}
\usepackage{mathbbol}
\usepackage[pagebackref]{hyperref}
\usepackage{color}
\usepackage{xcolor}
\usepackage{longtable}
\usepackage{lscape}
\usepackage{scalefnt}
\usepackage{setspace}
\usepackage{enumitem}
\usepackage{verbatim}
\usepackage[all]{xy}
\newcommand{\rar}{\rightarrow}
\newcommand{\lar}{\longrightarrow}
\newcommand{\llar}{-\kern-5pt-\kern-5pt\longrightarrow}
\newcommand{\surjects}{\twoheadrightarrow}

\newtheorem{Theorem}{Theorem}[section]
\newtheorem{Lemma}[Theorem]{Lemma}
\newtheorem{Corollary}[Theorem]{Corollary}
\newtheorem{Proposition}[Theorem]{Proposition}
\newtheorem{Remark}[Theorem]{Remark}
\newtheorem{Example}[Theorem]{Example}
\newtheorem{Conjecture}[Theorem]{Conjecture}
\newtheorem{Definition}[Theorem]{Definition}

\def\sqr#1#2{{\vcenter{\hrule height.#2pt
        \hbox{\vrule width.#2pt height#1pt \kern#1pt
            \vrule width.#2pt}
        \hrule height.#2pt}}}
\def\phi{\varphi}
\def\demo{\noindent{\bf Proof. }}
\def\square{\mathchoice\sqr64\sqr64\sqr{4}3\sqr{3}3}
\def\qed{\hspace*{\fill} $\square$}

\DeclareMathOperator{\depth}{depth}

\def\xx{{\bf x}}
\def\yy{{\bf y}}
\def\zz{{\bf z}}

\def\gg{{\bf g}}

\def\fm{{\mathfrak m}}
\def\fn{{\mathfrak n}}

\def\depth{{\rm depth}\,}

\def\restr{{\kern-1pt\restriction\kern-1pt}}

\def\R{{\mathbb R}}

\def\pp{{\mathbb P}}

\begin{document}
\begin{center}
\Large{\sc De Jonqui\`eres transformations in arbitrary dimension.
}\\
\large{\sc An ideal theoretic view}

	\vspace{0.3in}
	{\large\sc Zaqueu Ramos}\footnote{Partially supported by a CNPq grant (305860/2019-4))} 
	\quad\quad
	{\large\sc Aron  Simis}\footnote{Partially supported by a CNPq grant (301131/2019-8).}

\end{center}


\medskip


\begin{abstract}

A generalization of the plane de Jonqui\`eres transformation to arbitrary dimension is studied, with an eye for the  ideal theoretic side.
In particular, one considers structural properties of the corresponding base ideal and of its defining relations. Useful throughout is the idea of downgraded sequences of forms,  a tool considered in many sources for the rounding-up of ideals of defining relations.
The emphasis here is on the case where the supporting Cremona transformation of the de Jonqui\`eres transformation is the identity map. In this case we establish aspects of the homological behavior of the graph of the transformation.
	
\end{abstract}

\section*{Introduction}

Plane de  Jonqui\`eres transformations are  classically well-known.
Beside their simplicity amid other plane Cremona transformations, they constitute the essential tool in Castelnuovo's version of N\"other's theorem on the generation of the plane Cremona group.
The idea of introducing similar transformations in arbitrary dimension was probably present in the minds of many experts.
An explicit description of such an analogue appears in \cite{PanStellar} under the designation of {\em transformations stellaires} that reveals its geometric nature in terms of  {\em \'etoiles} of lines through a point.
For a group theoretic approach to these transformations in $\pp^3$ see 
 \cite{PanSim} and, for an elimination approach, see \cite{HS2}.
 
 The core of this generalization to arbitrary dimension relies on the existence of a given Cremona transformation in one lower dimension that acts as its {\em supporting Cremona} transformation.
 In the plane case this Cremona transformation is the identity map, so it becomes nearly invisible.
 However, already in $\pp^3$ the theory gets sufficiently entangled in some of its facets, specially if one takes the algebraic side of the transformations -- in the sense of ideal theory, not group theory.
In this work we develop ideal theoretic aspects of these transformations as birational maps with a discussion of the structure of the corresponding graphs (blowups).
Some of these aspects have been trailed in \cite{HS_Cremona} in the plane case, in quite some detail.
The parallel case of identity support in arbitrary dimension is dealt here with equally exhausting detail.
Our approach is rendered by applying the idea of a certain sequence of Sylvester like forms, here called {\em fully downgraded sequence}.
This idea is present in many places and has an elimination nature pattern (see, e.g., \cite{Cox}, \cite{Syl1}, \cite{BenDandr}, \cite{HS_Cremona},\cite{SiTo}, \cite{BuSiTo_Sylvester}). 

The reason to develop the ideal theoretic side of the geometric content of the subject, in addition to establishing a potential new research topic, is to stress a few recent techniques, such as the employ of downgraded sequences as mentioned above and the algebraic version of the fibers of a rational map by Eisenbud and Ulrich (\cite{EiUl}, also \cite{KPU}).

We give a brief description of the two sections.

In the first section we introduce the idea of two confluent transformations, as the socle of the theory. The geometric background has been established in \cite{PanStellar}, whereas here we develop the foundations of the algebraic side, such as a more precise facet of the corresponding base ideals of the transformations.
We next expand on the notion of a downgraded sequence, a tool that will be used both in the birational theory of two confluent rational map, leading to the notion of a generalized de Jonqui\`eres transformation $\mathfrak{J}$, and in the details of the defining relations of the graph of $\mathfrak{J}$ (blowup).
The rest of the section is devoted to the basics of the thus called generalized de Jonqui\`eres transformation $\mathfrak{J}$ with support  a Cremona transformation $\mathfrak{F}$. We show that this is a birational map whose inverse is again a map of the same sort, whose  support is the inverse of $\mathfrak{F}$.

The second section is devoted to the case where the supporting Cremona transformation of a generalized de Jonqui\`eres transformation $\mathfrak{J}$ is the identity.
We give the minimal free resolution of the  base ideal $I$ of $\mathfrak{J}$. 
In particular, having the details of the presentation ideal of the symmetric algebra of $I$, we head on to the structure of the presentation ideal of the Rees algebra of $I$, corresponding to the graph of $\mathfrak{J}$.
First, drawing again upon the fully downgraded sequences we give the explicit base ideal of the inverse of $\mathfrak{J}$ in terms of a certain defining equation of the graph of $\mathfrak{J}$. As a consequence we obtain that the inverse is defined in the same degree as $\mathfrak{J}$.
Next, as the main result of the section, we obtain a minimal set of generators of the presentation ideal of the graph in terms of the elements of a fully downgraded sequence.
As an application, we show that the Rees algebra of $I$ is an almost Cohen--Macaulay domain and conjecture that it is Cohen--Macaulay if and only if $d\leq n+1$, where $d\geq 2$ is the defining degree of $\mathfrak{J}$ and $n+1$ is the ambient number of variables.

\section{Main notions}

\subsection{Confluent transformations}

Let $\pp^n=\pp^n_k$ and $\pp^m=\pp^m_k$ stand for projective spaces over an infinite  field $k$.
A rational map $F:\pp^n\dasharrow \pp^m$ is defined by $m+1$ forms $f_1,\ldots,f_{m+1}$ of same degree in the source polynomial ring $R=k[x_1,\ldots,x_{n+1}]$.
In an  imprecise way, we write $F=(f_1:\cdots :f_{m+1})$, where the colon signals for homogeneous coordinates in $\pp^m_{k(x_1,\ldots,x_{m+1})}$. The ideal of $R$ generated by the defining forms is called the {\em base ideal} of the map.
The map $F$ is said to be birational onto its image if there is a rational map $\pp^m\dasharrow \pp^n$ whose restriction to the image of $F$ is an inverse. Birational maps of $\pp^n$ to itself are called Cremona transformations.

Let $\mathfrak{J}:\pp^n\dasharrow \pp^n$  be a rational map. 
As usual, for the algebraic discussion of such maps, we will distinguish $\pp^n_x={\rm Proj}(k[x_1,\ldots,x_{n+1}])$ (source) and $\pp^n_y={\rm Proj}(k[y_1,\ldots,y_{n+1}])$ (target).

Given a rational map $\mathfrak{F}:\pp^{n-1}\dasharrow \pp^{n-1}$, where,  similarly, $\pp^{n-1}_{x'}={\rm Proj}(k[x_1,\ldots,x_{n}])$ (source) and $\pp^{n-1}_{y'}={\rm Proj}(k[y_1,\ldots,y_{n}])$ (target), 
let $\pi_x: \pp^n_x\surjects \pp^{n-1}_{x'}$  and $\pi_y: \pp^n_y\surjects \pp^{n-1}_{y'}$ denote the natural coordinate projections, themselves rational maps.

\begin{Definition}\rm
	We say that $\mathfrak{J}$ and $\mathfrak{F}$ are {\em confluent} if
	\begin{equation}\label{basic_assumption_deJonq}
		\mathfrak{F}\circ \pi_x=\pi_y\circ \mathfrak{J}
	\end{equation}
	as rational maps from $\pp^n_x$ to $\pp^{n-1}_y$.
\end{Definition}
Since we deal quite a bit with inverses, we will obviously allow for a similar notion by reversing variables.


\begin{Lemma}\label{Jonq_basic_form}
	Let $\{g_1,\ldots, g_n \}\subset k[x_1,\ldots,x_{n}]$ be forms of the same degree with no proper common factor and let
	$\mathfrak{J}$ and $\mathfrak{F}$ be confluent rational maps as above.  Then $\mathfrak{F}=(g_1:\cdots : g_n)$  if and only if $\mathfrak{J}=(fg_1 :\cdots : fg_n:g)$ for suitable forms $f,g\in k[x_1,\ldots,x_{n+1}]$ satisfying $\gcd(f,g)=1$.
\end{Lemma}
\demo 
Let $\mathfrak{F}=(g_1 :\cdots : g_n)$ and set $\mathfrak{J}=(h_1 :\cdots : h_n: h_{n+1})$, with $\gcd(h_1,\ldots,h_n,h_{n+1})=1$.
Then
$$\mathfrak{F}\circ \pi_x=(g_1(x_1,\ldots,x_n)  :\cdots :  g_n(x_1,\ldots,x_n))=(g_1  :\cdots :  g_n)$$
and
$$\pi_y\circ \mathfrak{J}=(y_1(h_1,\ldots ,h_n,h_{n+1}):\cdots :y_n(h_1,\ldots ,h_n,h_{n+1}))=(h_1:\cdots :h_n).
$$
Therefore, $(h_1  :\cdots : h_n)=(g_1  :\cdots : g_n)$ define the same rational map, hence  the corresponding vectors are proportional.
This means that there exist nonzero forms $p,q\in k[x_1,\ldots,x_{n+1}]$ with the right degrees such that $ph_i=qg_i$ for $1\leq i\leq n$.
But, since $\gcd(g_1,\ldots,g_n)=1$, then $p=1$.
Thus,
$$(h_1  :\cdots : h_n ; h_{n+1})=(qg_1  :\cdots : qg_n :h_{n+1}).$$
Now, take $f:=q$ and $g:=h_{n+1}$.
Finally, since $f$ divides every $h_i \,\, (1\leq i\leq n)$ and $\gcd(h_1,\ldots,h_n,h_{n+1})=1$, then $\gcd(f,g)=1$.

Conversely, let $\mathfrak{J}=(fg_1 :\cdots :  fg_n ; g)$ with $\gcd(f,g)=1$.
Then
\begin{eqnarray*}
	\mathfrak{F}\circ \pi_x\kern-12pt&=&\kern-10pt\pi_y\circ \mathfrak{J}=(y_1(fg_1 , \ldots , fg_n, g)  :\cdots :  y_n(fg_1 , \ldots , fg_n, g))=(fg_1  :\cdots : g_n)\\
	&=& \kern-10pt(g_1  :\cdots : g_n)
\end{eqnarray*}\
as a rational map.
Therefore, $\mathfrak{F}$ is defined by $g_1,\ldots,g_n$.
\qed

\smallskip

A preliminary tool related to the above environment is the following:

\begin{Lemma}\label{arealgind} With the above notation, suppose that $g_1,\ldots,g_n$ are algebraically indeependent over $k$ and that $\gcd(f,g)=1$. Then  $fg_1,\ldots,fg_n,g$ are algebraically independent over $k$.	
\end{Lemma}
\demo
The argument is an  elementary calculation using the fact that $\gcd(f,g)=1$.
Namely, suppose   $F(y_1,\ldots,y_{n+1})\in k[y_1,\ldots,y_n]$ is a nonzero homogeneous polynomial such that $F(fg_1,\ldots,fg_n,g)=0.$ 
Write $F$ as a polynomial in $y_{n+1}$ with coefficients in $k[y_1,\ldots,y_n]$, say,  $F=p_{m}y_{n+1}^{m}+\cdots+p_1y_{n+1}+p_0$, with $p_i\in k[y_1,\ldots,y_n]$ homogeneous polynomials of degree $\deg F-i$ with $p_m,p_0\neq 0.$  
Thus, the vanishing of $F(fg_1,\ldots,fg_n,g)$ translates into
$$f^{\deg F-m}[p_m(g_1,\ldots,g_n)g^m+\cdots+p_i(g_1,\ldots,g_n)f^{m-i}g^i+\cdots+p_{0}(g_1,\ldots,g_n)f^{m}]=0,$$
or, equivalently,
$$p_m(g_1,\ldots,g_n)g^m+\cdots+p_i(g_1,\ldots,g_n)f^{m-i}g^i+\cdots+p_{0}(g_1,\ldots,g_n)f^{m}=0.$$
Hence, $[p_m(g_1,\ldots,g_n)\,\cdots\, p_0(g_1,\ldots,g_n)]^t$ is a syzygy of the ideal power $(f,g)^m$, where  $p_m(g_1,\ldots,g_n)\neq 0$ and $p_0(g_1,\ldots,g_n)\neq 0$ since $g_1,\ldots,g_n$ are algebraically independent over $k$.
But the syzygy matrix of $(f,g)^m$ is well-known to be the following  $(m+1)\times m$ matrix
	$$
	\arraycolsep=3pt 
	\medmuskip = 4mu
	\left(\begin{array}{cccccccccc}-g&\\
	f&-g&\\
	&f&\ddots\\
	&&\ddots\\
&&&-g\\
&&&f\end{array}\right).$$ 
Thus, $p_m(g_1,\ldots,g_n)=pf^m$ and $p_0(g_1,\ldots,g_n)=qg^m$ for certain $p,q\in k[x_1,\ldots,x_{n+1}].$ Since $p_m(g_1,\ldots,g_n)\neq 0$ and $p_0(g_1,\ldots,g_n)\neq 0$ and neither depends depend on $x_{n+1}$, while either $f$ or $g$ does, we arrive at  a contradiction. 
\qed

\subsection{Fully downgraded sequences}\label{Sylvester_basics}

The idea of introducing certain downgraded forms to get insight into the structure of generators of the blowup equations has appeared in many places (see \cite{Cox}, \cite{Syl1}, \cite{BenDandr}, \cite{HS_Cremona},\cite{SiTo}, \cite{BuSiTo_Sylvester} for a few recent ones).

Introduce additional notation, as follows: $R':=k[\xx']=k[x_1,\ldots,x_n]$ and $R:=k[\xx]=k[x_1,\ldots,x_n,x_{n+1}]$.
As in the previous section, $g_1,\ldots,g_n\in R'$ are forms of the same degree $d'$ with no proper common factor and $f,g\in R:=k[\xx]=k[x_1,\ldots,x_n,x_{n+1}]$ are additional forms of degrees $d-d'$ and $d> d',$ respectively,  such that $\gcd(f,g)=1$ and, moreover,  one  at least involves $x_{n+1}$ effectively. 
Set, in addition, ${\bf g}:=\{g_1,\ldots,g_n\}$ and $f{\bf g}:=\{fg_1,\ldots,fg_n\}.$

In the sequel we will emphasize $f,g$ as polynomials in $x_{n+1}$ with coefficients in $k[\xx']$, in particular, the respective $x_{n+1}$-degree will come up quite often.

Let $I\subset R$ denote the base ideal generated by  $fg_1,\ldots,fg_n,g.$
We are interested in the Rees algebra $\mathcal{R}(I)=R[fg_1t,\ldots,fg_nt,gt]\subset R[t]$ and, more precisely, in its presentation ideal defined as the kernel $\mathcal{I}$ of the homomorphism of bigraded $k$-algebras
$$S:=R[y_1,\ldots,y_{n+1}]\surjects\mathcal{R}(I),\quad x_i\mapsto x_i,\,y_j\mapsto fg_jt \,(1\leq j\leq n),\,y_{n+1}\mapsto gt.$$
Consider the extended ideal $J:=(g_1,\ldots,g_n)R\subset R$.
Since the ideals $J,(fg_1,\ldots,fg_n)\subset R$ have the same syzygies, there is an isomorphism of the respective Rees algebras. On the other hand, $\mathcal{R}_R(fg_1,\ldots,fg_n)\subset \mathcal{R}_R(I)$ induces an inclusion of the respective presentation ideals.

Assume, moreover, that $g_1,\ldots,g_n$ define a Cremona map $\mathfrak{F}$ of $\pp^{n-1}.$ 

Let ${\bf h}=h_1,\ldots,h_n\in k[y_1,\ldots,y_n]$  be forms of degree $d'$ defining the   inverse  map $\mathfrak{F}^{-1}$ and let  $\Delta$ be the corresponding source inversion factor, i.e., the uniquely defined form of $R'$  of degree $dd'-1$ such that
\begin{equation}\label{source_inversion_factor}
	h_i({\bf g})=\Delta\cdot x_i \quad \quad (i=1,\ldots,n)
\end{equation}


Now let $\mathfrak{z}$ be a syzygy of $f{\bf g},g$ of degree $e$ whose last coordinate is nonzero -- e.g., any Koszul syzygy involving $g$. Say, $\mathfrak{z}$  has standard degree $e$. Denote by $\delta$ the largest integer such that $(x_1,\ldots,x_n)^\delta$ contains every entry of  $\mathfrak{z}$ -- clearly, $\delta\leq e$.
Then $F_1:=[y_1\,\cdots\, y_{n+1}]\cdot\mathfrak{z}\in \mathcal{I}$  is a form of bidegree $(e,1)$ and,  moreover,  $F_1\in (x_1,\ldots,x_n)^{\delta}S$.  Thus,  write 
$F_1=F_{1,1}x_1+\cdots+F_{1,n}x_n,$
where $F_{1,i}$ are bihomogeneous polynomials in $(x_1,\ldots,x_n)^{\delta-1}S$   of bidegree $(e-1,1).$ 
Define

$$F_2=F_{1,1}h_1+\cdots+ F_{1,n}h_n.$$
In particular, $F_2$ has bidegree  $(e-1,d'+1)$. Moreover:
\begin{eqnarray}
	F_2(\xx,f{\bf g},g)&=&F_{1,1}(\xx,f{\bf g},g)h_1(f{\bf g})+\cdots+ F_{1,n}(\xx,f{\bf g},g)h_n(f{\bf g})\nonumber\\
	&=&f^{d'}\Delta[F_{1,1}(\xx,f{\bf g},g)x_1+\cdots+ F_{1,n}(\xx,f{\bf g},g)x_n] \quad\quad\quad\quad\quad(\mbox{by}\, \eqref{source_inversion_factor})\nonumber\\
	&=&f^{d'}\Delta\cdot F_1(\xx,f{\bf g},g)=0\quad\quad\quad\quad\quad\quad \quad\quad\quad\quad\quad\quad\quad\quad\quad(\mbox{because}\, F_1\in\mathcal{I})\nonumber
\end{eqnarray}
Thus, $F_2\in \mathcal{I}_{(e-1,d'+1)}.$

Note that $F_2$ is obtained from $F_1$ via replacing the template of $x_1,\ldots,x_n$ by the template of $h_1,\ldots,h_n$.
Moreover,  $h_1,\ldots,h_n$ are algebraically independent, as they define (the inverse of) a Cremona map.
Therefore,  $F_2$ is a nonzero $y_{n+1}$-monoid.



Again, $F_{1,i}\in (x_1,\ldots,x_n)^{\delta-1}S$, so we can write 
$F_2=F_{2,1}x_1+\cdots+F_{2,n}x_n$
where $F_{2,i}$ are bihomogeneous polynomials in $(x_1,\ldots,x_n)^{\delta-2}S$   of bidegree $(e-2,d'+1).$ 
Define $$F_3:=F_{2,1}h_1+\cdots+F_{2,n}h_n.$$ As above, $F_3$ is a form of bidegree $(e-2,2d'+1)$ and
$$ F_3(\xx,f{\bf g},g):=f^{d'}\Delta\cdot F_2(\xx,f{\bf g},g)=0.$$
Thus, $F_3\in  \mathcal{I}_{(e-2,2d'+1)}.$  We also have, similarly for $F_2,$ that $F_3$ is a a nonzero $y_{n+1}$-monoid. 
Iterating this process eventually gives a nonzero form $F_{\delta}\in \mathcal{I}$ of bidegree $(e-\delta,\delta d'+1).$ 

\begin{Definition}\label{Sylvester_seq}\rm
The set $\{F_1,F_2,\ldots, F_{\delta}\}$ is a {\em  fully downgraded sequence} associated to the syzygy $\mathfrak{z}$ of $I=(f\gg,g)$, where $\gg$ defines a Cremona map of $\pp^{n-1}$.
\end{Definition}
This sort of sequence, obtained by suitably downgrading the degrees of the $\xx$-terms in benefit of upgrading the ones of the $\yy$-terms, has been used in several occasions (see, e.g., \cite{HS2}).

\begin{Corollary}\label{formbidegree1D}
Keeping the above notation, suppose that $\gg:=\{g_1,\ldots,g_n\}$ define a Cremona map of $\pp^n$  and that $f,g$ are $x_{n+1}$-monoids with no proper common factor, one of which at least involves $x_{n+1}$ effectively.
	Then the ideal $I=(f\gg,g)$ admits a nonzero Rees equation $F$ of bidegree $(1,D)$, for suitable $D\geq 1$, depending effectively  on $x_{n+1}$.
\end{Corollary}
\demo Since $f,g$ are $x_{n+1}$-monoids then there is a syzygy $\mathfrak{z}$ of degree, say,  $\delta+1$ such that  $\delta$ is the largest integer for which the ideal $(x_1,\ldots,x_n)^\delta$ contains the entries of  $\mathfrak{z}.$  Hence, the term $F_{\delta}$ of the Sylvester sequence associated to $\mathfrak{z}$ has bidegree $(\delta+1-\delta,\delta d'+1)=(1,D)$, with $D=\delta d'+1$.
Now, $F_{\delta}$ depends effectively on $x_{n+1}$,  as otherwise,  we could produce a  term  $0\neq F_{\delta +1}\in \mathcal{I}$ of bidegree $(0, \delta d'+2),$ which is an absurd  by Lemma~\ref{arealgind}.
\qed

\subsection{Generalized de Jonqui\`eres transformations}

Keep the notation  of Lemma~\ref{Jonq_basic_form},.
We give an algebraic proof of the basic structural result of \cite{PanStellar}.

\begin{Proposition}\label{generalized_deJonq_setup}
	Suppose that $\mathfrak{F}$ and $\mathfrak{J}$ are confluent, where $\mathfrak{F}=(g_1  :\cdots : g_n)$, with $g_1,\ldots,g_n$ having no proper common factor, and $\mathfrak{J}=(fg_1 :\cdots :fg_n :g)$, with $\gcd(f,g)=1$. 
	The following are equivalent:
	\begin{enumerate}
		\item[{\rm (i)}]  $\mathfrak{J}$ is birational.
		\item[{\rm (ii)}]  $\mathfrak{F}$ is birational  and, moreover,  $f,g$ are $x_{n+1}$-monoids  one of which at least involves $x_{n+1}$ effectively.
	\end{enumerate}
\end{Proposition}
\demo
(i) $\Rightarrow$ (ii) 
Let $\mathfrak{i}_n$ denote the identity map of $\pp^n$, regardless as to whether it is source or target -- the details of calculation will readily identify which is which.
Similarly, $\mathfrak{i}_{n-1}$ is the identity map of $\pp^{n-1}$.

By the present assumption, $\mathfrak{J}\circ \mathfrak{J}^{-1}=\mathfrak{J}^{-1}\circ \mathfrak{J}=\mathfrak{i}_n$.
The standing assumption implies
\begin{equation}\label{relation1}
	\mathfrak{F}\circ \pi_x\circ \mathfrak{J}^{-1}= \pi_y\circ \mathfrak{J}\circ \mathfrak{J}^{-1}=\pi_y\circ \mathfrak{i}_n=\pi_y.
\end{equation}
On the other hand, letting $L_y=V(y_{n+1})\subset \pp^n$, one has an isomorphism $\pi_y|_{L_y}: L_Y\simeq \pp^{n-1}$ induced by restriction of the coordinate projection $\pi_y$.
Let $(\pi_y|_{L_y})^{-1}$ denote its inverse.
Applying (\ref{relation1})  we find
$$(\mathfrak{F}\circ \pi_x\circ \mathfrak{J}^{-1})|_{L_y}=\pi_y|_L,
$$
which further implies
$$\mathfrak{F}\circ(\pi_x\circ \mathfrak{J}^{-1}|_{L_y}\circ (\pi_y|_{L_y})^{-1})=\pi_y|_{L_y}\circ (\pi_y|_{L_y})^{-1}=\mathfrak{i}_{n-1}. 
$$
Therefore, $\mathfrak{F}$ has an inverse.

For the additional assertion we argue that:

{\sc Claim 1.}   If $\deg_{x_{n+1}} (f)=\deg_{x_{n+1}} (g)=0$, the map $\mathfrak{J}$ is not dominant.

This is clear as the defining forms of $\mathfrak{J}$ do not involve $x_{n+1}$, triggering 
$$\dim k[fg_1,\ldots,fg_n,g]\leq n.$$

{\sc Claim 2.} $\deg_{x_{n+1}} (f)\leq 1$ and $\deg_{x_{n+1}} (g)\leq 1$.

Set $\deg_{x_{n+1}} (f)=r$ and $\deg_{x_{n+1}} (g)=s$.

Let $\boldsymbol\alpha=(\alpha_1:\cdots:\alpha_n:\alpha_{n+1})\in \pp^n_x$ be a general point  in the source and let $\boldsymbol\alpha'=(\alpha_1:\cdots:\alpha_n)\in \pp^{n-1}_{x'}$.
Since gcd$(f,g)=1$ and $\boldsymbol\alpha$ is general then gcd$(f(\boldsymbol\alpha',x_{n+1}),g(\boldsymbol\alpha',x_{n+1}))=1.$ Hence, $f(\boldsymbol\alpha',x_{n+1}),g(\boldsymbol\alpha',x_{n+1})$ have no common root as polynomials in $x_{n+1}$. In particular,  the polynomial $$p=p(x_{n+1}):=g(\boldsymbol\alpha)f(\boldsymbol\alpha',x_{n+1})-f(\boldsymbol\alpha)g(\boldsymbol\alpha',x_{n+1})\in k[x_{n+1}]$$
has no common root with either $f(\boldsymbol\alpha',x_{n+1})$ or $g(\boldsymbol\alpha',x_{n+1}).$

Moreover, one has:

{\sc Claim 3.} $\deg p= \max\{r,s\}$.

Clearly, $\deg p\leq \max\{r,s\}$. If equality fails one must have $r=s$.
Let $\mathfrak{f}, \mathfrak{g}\in k[x']$ respectively denote the $r$th coefficients of $f,g$ as polynomials in $x_{n+1}$ over $k[\xx']$.
Then the $r$th coefficient of $p$ in $k[x_{n+1}]$ is $g(\boldsymbol\alpha)\mathfrak{f}-f(\boldsymbol\alpha)\mathfrak{g}$.
Now, $g\mathfrak{f}-f\mathfrak{g}\neq 0$ since $\gcd(f,g)=1$ entails that any syzygy of $\{f,g\}$ is a multiple of its Koszul syzygy, while $\mathfrak{f}, \mathfrak{g}$ do not depend on $x_{n+1}$.
Since $\boldsymbol\alpha$ is a general point then also $g(\boldsymbol\alpha)\mathfrak{f}-f(\boldsymbol\alpha)\mathfrak{g}\neq 0$, as was claimed.

To proceed, it will suffice to show that $\deg p=1$.

For this, let $\boldsymbol\beta=(\beta_1:\cdots :\beta_{n+1})\in \pp^n_y$ be the point in the target which is is the image by $\mathfrak{J}$ of the point $\boldsymbol\alpha$, a  well-defined value as  $\boldsymbol\alpha$ is general.
According to \cite[Definition 3.3]{KPU}, the fiber (inverse image of $\boldsymbol\beta$ by $\mathfrak{J}$) is defined by the ideal
$$\mathbb{f}(\alpha):=I_2\left(\begin{array}{cccc}fg_1&\cdots&fg_n&g\\
\beta_1&\cdots&\beta_n&\beta_{n+1}\end{array}\right):I^{\infty}.$$

Since  $\mathfrak{J}$ is assumed to be birational, this fiber is reduced to a single proper point of multiplicity one. Thus, $\mathbb{f}(\alpha)$ is one-dimensional unmixed, hence regular.
Therefore, it coincides with the prime ideal of the point $\alpha$. Since this prime ideal is generated by ($n$ independent linear forms among) the $2$-minors of $I_2\left(\begin{array}{cccc}x_1 &\cdots&x_n&x_{n+1}\\
	\alpha_1&\cdots&\alpha_n&\alpha_{n+1}\end{array}\right)$, in particular one has $l:=\alpha_1x_{n+1}-\alpha_{n+1}x_1\in \mathbb{f}(\alpha)$.
Since $g\in I$,  then 
\begin{equation}\label{linear_form_inside_minors}
lg^t\in I_2\left(\begin{array}{cccc}fg_1&\cdots&fg_n&g\\
	\beta_1&\cdots&\beta_n&\beta_{n+1}\end{array}\right),
\end{equation}
for certain $t\geq 1$.

Now recall that $\beta_i=f(\boldsymbol\alpha)g_i(\boldsymbol\alpha')$, for $1\leq i\leq n$, while $\beta_{n+1}=g(\boldsymbol\alpha)$.
Thus, by evaluating $\xx'\mapsto \boldsymbol\alpha'$ all $2$-minors of the right side of (\ref{linear_form_inside_minors}) involving the first $n$ columns vanish, while the ones involving the last column are all multiple of $p=p(x_{n+1})$.
Therefore, the result of so evaluating (\ref{linear_form_inside_minors}) is that $l g(\boldsymbol\alpha',x_{n+1})^t$ is a multiple of $ p(x_{n+1})$  in the polynomial ring $k[x_{n+1}]$.
Since $p(x_{n+1})$ has no common root with $g(\boldsymbol\alpha',x_{n+1}),$ then $\deg p(x_{n+1})=1,$  as was to be shown.

\smallskip

(ii) $\Rightarrow$ (i) Let $B$ an $n\times (n-1)$ submatrix of the Jacobian dual matrix of $g_1,\ldots,g_n$ having rank $n-1.$ By Corollary~\ref{formbidegree1D} there is a form $F$ of bidegree $(1,D)$ in the defining ideal of the Rees algebra $\mathcal{R}(I)$ such that $\partial F_d/\partial x_{n+1}\neq 0.$  Hence, the matrix

$$\left(\begin{array}{cc|cccc}
	&&\partial F/\partial x_1&\\
	B&&\vdots\\
	&&\partial F/\partial x_n\\
	\hline &&&\\ [-7pt]
	\boldsymbol0&& \partial F/\partial x_{n+1}
\end{array}\right)$$
is an $(n+1)\times n$ submatrix of the jacobian dual matriz of $fg_1,\ldots,fg_n,g$ having rank $n.$ Hence, by \cite[Theorem 2.18]{AHA}, $fg_1,\ldots,fg_n,g$ define a Cremona map of $\pp^{n}.$
\qed

\begin{Definition}\label{deJonq_gen_def}\rm
	Under the equivalent conditions of Proposition~\ref{generalized_deJonq_setup}, the Cremona transformation $\mathfrak{J}$ is called a {\em generalized de Jonqui\`eres transformation}, with {\em Cremona support} $\mathfrak{F}$.
	In particular, every generalized de Jonqui\`eres transformation is confluent with its Cremona support.
\end{Definition}

\begin{Proposition} \label{generalized_inverse}
	The inverse transformation of a generalized de Jonqui\`eres transformation $\mathfrak{J}$ with Cremona support $\mathfrak{F}$ is a generalized de Jonqui\`eres transformation with Cremona support $\mathfrak{F}^{-1}$
\end{Proposition}
\demo Multiply both members of the basic equation (\ref{basic_assumption_deJonq}) by $\mathfrak{F}^{-1}$ on the left and by $\mathfrak{J}^{-1}$ on the right, to get $\mathfrak{F}^{-1}\circ \pi_y=\pi_x\circ \mathfrak{J}^{-1}$.
Therefore, $\mathfrak{J}^{-1}$ and $\mathfrak{F}^{-1}$ are confluent.
Then, by Lemma~\ref{Jonq_basic_form},  one may set $$\mathfrak{F}^{-1}=(g'_1  :\cdots : g'_n) \quad \text{\rm and}\quad \mathfrak{J}^{-1}=(f'g'_1  :\cdots : f'g'_n : g'),$$ 
for certain forms $g'_1,\ldots,g'_n, f',g'$ in $ k[y_1,\ldots,y_n ,y_{n+1}]$, with $\gcd(f',g')=1$. 
Thus, since $\mathfrak{J}^{-1}$ is birational,  Proposition~\ref{generalized_deJonq_setup} implies that 
 $f',g'$ are $y_{n+1}$-monoids with at least one of them effectively involving $y_{n+1}$.
 This shows that $\mathfrak{J}^{-1}$ is a generalized de Jonqui\`eres transformation with Cremona support $\mathfrak{F}^{-1}$.
\qed

\section{The case of  identity support}\label{deJonq_identity_support}

Changing a bit our original notation, let $\mathfrak{F}$ denote a generalized de Jonqui\`eres transformation of $\pp^n$, with Cremona support the identity $(x_1:\ldots :x_n)$.
Let $I=(fx_1,\ldots,fx_n,g)\subset R=k[x_1,\ldots,x_n,x_{n+1}]$ stand for the corresponding base ideal.

\subsection{The minimal free resolution of the base ideal}

A free (not necessarily minimal) resolution of the base ideal $I$ is available for arbitrary $n$ and arbitrary Cremona support by a straightforward adaptation of the result in \cite[Proposition 2.3]{HS2}.
As a matter of fact, this resolution is kosher regardless of the monoid condition.
We give a proof in the case of identity support for the reader's convenience.
We assume throughout  that $g\in (x_1,\ldots,x_n)R$, where $R=k[x_1,\ldots,x_{n},x_{n+1}]$, i.e., the $x_{n+1}$-coefficient of $g$ is not a scalar.

\begin{Proposition}\label{resolution_of_deJonq}
Let $f,g\in R=k[x_1,\ldots,x_{n},x_{n+1}]$ be relatively prime forms  of respective degrees $d-1$ and $d\geq 2$ and such that  $g\in (x_1,\ldots,x_n)R$. Set $I:=(x_1f,\ldots,x_{n}f,g)$. Then the minimal graded free resolution of $R/I$ is
	$$\small{0\to R(-(n+d-1))\to\cdots\to R(-(d+2))^{n\choose3}\to\begin{array}{c}R(-(d+1))^{n\choose2}\\\oplus\\R(-(2d-1))\end{array}\stackrel{\psi}\lar R(-d)^{n+1}\to R,
	}$$
	where
	\begin{enumerate}
		\item[{\rm (a)}] The syzygy matrix $\Psi$ has the form
		$$\psi=\left(\begin{array}{ccc|ccccccc}&&&-q_1\\
			&\mathcal{K}&&\vdots\\
			&&&-q_{n}\\
			\hline
			&&&\\ [-7pt]
			0&\cdots&0&f\end{array}\right)$$
		with $\mathcal{K}$ denoting the first Koszul matrix of  $x_1,\ldots,x_{n}$ and $g=q_1x_1+\cdots +q_nx_n$. 
		\item[{\rm (b)}]  The $n-2$ tail terms are the $n-2$ tail terms of the Kozul complex of $\{x_1,\ldots,x_n\}$ shifted by $d-1$.
	\end{enumerate}
\end{Proposition}
\demo Set $\fm=(x_1,\ldots,x_n)$. The minimal graded free resolution of $R/f\fm$ is the Koszul complex of the generators of $\fm$ shifted by $\deg f=d-1$:
\begin{equation}\label{koszul}
	\small{0\to R(-(n+d-1))\to\cdots 
		\to R(-(d+1))^{n\choose2}\stackrel{\psi}{\to} R(-d)^{n}\to R \to R/f\fm\to 0.}
\end{equation}
Now, $f\fm:g=(f)$. Shifting by $d$, we get the exact sequence
\begin{equation}\label{res_principal}
	0\to R(-(2d-1))\to R(-d)\to (R/f\fm:g)(-d) \to 0.
\end{equation}
On the other hand, multiplication by $g$ on $R$ induces the exact 
sequence $$0\to (R/f\fm:g)(-d) \stackrel{\cdot g}{\to}  R/I'\to R/I\to 0.$$
Extending to a map from \eqref{koszul}  to \eqref{res_principal},
the resulting mapping cone
$$\small{0\to R(-(n+d-1))\to\cdots\to R(-(d+2))^{n\choose3}\to\begin{array}{c}R(-(d+1))^{n\choose2}\\\oplus\\R(-(2d-1))\end{array}\stackrel{\psi}\lar R(-d)^{n+1}\to R}$$
is a graded free resolution of $R/I.$ 
Due to its shifts, this resolution  is actually minimal.
\qed

\begin{Corollary}\label{saturation_CM}
	With the notation of {\rm Proposition~\ref{resolution_of_deJonq}}, one has:
	\begin{enumerate}
		\item[{\rm (i)}] $I$ is saturated.
		\item[{\rm (ii)}] $(x_1,\ldots,x_n)$ is an associated prime of $R/I$.
		\item[{\rm (iii)}]  $R/I$ is Cohen--Macaulay if and only if $n=2$, i.e., when $I$ defines a plane de Jonqui\`eres transformation.
	\end{enumerate}
\end{Corollary}
\demo
(i) It follows from the resolution which gives $\depth R/I=1$.

(ii) Since $(x_1,\ldots,x_n)\subset I:f$, then $(x_1,\ldots,x_n)$ is contained in an associated prime of $R/I$.
If the containment were proper it would contradict (i).

(iii) If $R/I$ is Cohen--Macaulay, it follows from (ii) that $n=2$. Conversely, if $n=2$ then $I$ is the ideal of maximal minors of a $3\times 2$ matrix (\cite[Proposition 2.3]{HS_Cremona}).
\qed

\begin{Example}\rm (Plane de Jonqui\`eres transformation) \label{deJonq_plane}
The generalized de Jonqui\`eres transformation of dimension $2$, with identity support, is the classical {\em de Jonqui\`eres transformation} with homaloidal type $(d;d-1,1^{2d-2})$, for any given $d\geq 1$.
It is in fact the only generalized de Jonqui\`eres transformation for $n=2$.
This is because the only Cremona map of $\pp^1$ is the identity.

The base ideal $I\subset k[x,y,z]$ of a plane  de Jonqui\`eres transformation of degree $d\geq 2$  has the following properties:

$\bullet$ $k[x,y,z]/I$ is a Cohen--Macaulay ring of multiplicity  $d(d-1)+1.$

$\bullet$ $I$ is an ideal of linear type if and only if $d=2.$
\end{Example}

\subsection{The graph}

Recall from Subsection~\ref{deJonq_identity_support} that a generalized de Jonqui\`eres transformation has base ideal $I=(x_1f,\ldots,x_nf,g)$, where $\deg(g)=d\geq 2$ and $f,g$ are $x_{n+1}$-monoids satisfying an additional condition.
In particular,  $g\in (x_1,\ldots,x_n)R$, with $R=k[x_1,\ldots,x_{n+1}]$.

By Proposition~\ref{resolution_of_deJonq} (a), the syzygy matrix of $I$ has the following shape
$$\psi=\left(\begin{array}{ccc|ccccccc}&&&-q_1\\
	&\mathcal{K}&&\vdots\\
	&&&-q_{n}\\
	\hline
	&&&\\ [-7pt]
	0&\cdots&0&f\end{array}\right)$$
where $\mathcal{K}$ is the first Koszul matrix of  $x_1,\ldots,x_{n}$ and $g=q_1x_1+\cdots x_n$. 
The  forms 
\begin{equation}\label{symm_algebra_equations}
p_{i,j}:=x_jy_i-x_iy_j \,\,(1\leq i<j\leq n)\quad\mbox{and}\quad F:=f\,y_{n+1}-\displaystyle\sum_{i=1}^{n} q_i. y_i
\end{equation}
associated to the columns of $\psi$ belong to the presentation ideal $\mathcal{I}\subset R[y_1,\ldots,y_{n+1}]$ of the Rees algebra $\mathcal{R}_R(I)$ of $I$, obtained by mapping
$y_i\mapsto x_if$, for $1\leq i\leq n$, and $y_{n+1}\mapsto g$.

\medskip

\subsubsection{Main structure result}

A fully downgraded sequence is much simpler when the Cremona support is the identity.
 There is an obvious syzygy to play the role of $\mathfrak{z}$ in the notation of Subsection~\ref{Sylvester_basics}, and that is $F$ as in (\ref{symm_algebra_equations}); it has standard degree $\deg(f)=\deg(g)-1=d-1$.

Let then $\{F_0=F, F_1,\ldots, F_{d-2}\}\subset \mathcal{I}$ stand for a fully downgraded sequence, where for convenience we have shifted the indices by one to start at zero.
Recall that, by construction, in this case $F_{d-2}$ has bidegree $(1,d-1)$, since $d-1$ is the standard degree of the right most column of $\psi$.
We state the main highlights about these forms in the next remark.

\begin{Remark}\label{destaques}\rm
	The forms $F_i$, $0\leq i \leq d-2, $ have the following properties :
	\begin{enumerate}
		\item[(a)]  $F_i$ is a strict $y_{n+1}$-monoid in $k[x_1,\ldots,x_{n+1}, y_1, \ldots, y_{n+1}]$.
		\item[(b)]  If $d\geq 3$ then $F_i\in (x_1,\ldots,x_n)^{d-(i+1)}R[y_1,\ldots,y_{n+1}]$, for $0\leq i\leq d-2.$ 
		\item[(c)] 
		$x_jF_{i}-y_jF_{i-1}=\sum (x_jy_k-x_ky_j)F_{i-1,k}\quad (j=1,\ldots,n)$, for $1\leq i\leq d-2$.
	\end{enumerate}
\end{Remark}
Only (a) needs a bit of an argument. Namely, for any chosen decomposition $F_0=F_{0,1}x_1+\cdots+F_{0,n}x_n$ as displayed in Subsection~\ref{Sylvester_basics}, there is a unique $1\leq i\leq n$ such that $F_{0,i}=f\,y_{n+1}/x_i$.
This nonzero form $F_{0,i}$ is then the only entry of the corresponding ``content matrix''  having a term in $y_{n+1}$.
Therefore, the determinant has degree one in $y_{n+1}$.
The argument for the subsequent $F_i$'s is the same.

Using $F_{d-2}$ we can prove the following naturally expected property:

\begin{Proposition}\label{inverse_of_deJonq}
	Let $\mathfrak{J}$ denote a generalized de Jonqui\`eres map of degree $d\geq 2$ with identity support.
	Then the inverse map is defined by the forms
	\begin{equation*}
		\frac{\partial F_{d-2}}{\partial x_{n+1}}\,y_1,\ldots,\, \frac{\partial F_{d-2}}{\partial x_{n+1}}\,y_n, \,\sum_{i=1}^n \frac{\partial F_{d-2}}{\partial x_i}\, y_i.
	\end{equation*}
	In particular, the inverse map is a generalized de Jonqui\`eres map of degree $d$ with identity support.
\end{Proposition}	
\demo  An $(n+1)\times n$ submatrix of the Jacobian dual matrix of $\mathfrak{F}$ is 
\begin{equation}\label{syzygy_of_deJonq_dual}
	\Psi=\left(\begin{array}{ccc|ccccccc}&&&\partial F_{d-2}/\partial x_1\\
		&\mathcal{K}(\yy')&&\vdots\\
		&&&\partial F_{d-2}/\partial x_n\\[4pt]
		\hline
		&&&\\ [-7pt]
		0&\cdots&0&\partial F_{d-2}/\partial x_{n+1}
	\end{array}\right),
\end{equation}
where $\mathcal{K}(\yy')$ denotes the first Koszul matrix of the variables $\yy'=y_1,\ldots,y_n$.
Set 
$$\mathfrak{f}:=\frac{\partial F_{d-2}}{\partial x_{n+1}} \quad\mbox{and}\quad  \mathfrak{g}:= \sum_{i=1}^n \frac{\partial F_{d-2}}{\partial x_i}\, y_i.$$
 By the principle of the Jacobian dual method (\cite[Theorem 2.18 (i)]{AHA}), the ordered signed $n$-minors of $\Psi$ define  the inverse map  $\mathfrak{J}^{-1}.$ The latter are $y_n^{n-1}\mathfrak{f}y_1,\ldots,y_n^{n-1}\mathfrak{f}y_n,y_n^{n-1}\mathfrak{g}.$ Hence, $\mathfrak{f}y_1,\ldots,\mathfrak{f}y_n,\mathfrak{g}$ also define $\mathfrak{J}^{-1}.$  
 To conclude, we prove the following

{\sc Claim.}  $\mathfrak{f}$ and $\mathfrak{g}$ have no proper common factor.

Set $\mathfrak{h}={\rm gcd}(f,g).$ Then, for $\mathfrak{f}'=\mathfrak{f}/\mathfrak{h}$ and $\mathfrak{g}'=\mathfrak{g}/\mathfrak{h},$ the forms $\mathfrak{f}'y_1,\ldots,\mathfrak{f}'y_n,\mathfrak{g}'$ define also a representative to $\mathfrak{J}^{-1}.$ Say that $d':=\deg\mathfrak{f}'y_i=\deg \mathfrak{g}'.$ By definition, $d'\leq d.$ Let $\mathcal{I}'$ be presentation ideal of the Rees algebra of the ideal $(\mathfrak{f}'y_1,\ldots,\mathfrak{f}'y_n,\mathfrak{g}')\subset k[y_1,\ldots,y_n,y_{n+1}]$. 
By the above discussion, there is a biform $F'_{d'-2}\in\mathcal{I}'_{(d'-1,1)}$  effectively depending on $y_{n+1}.$ But, $\mathcal{I}'=\mathcal{I}.$ Hence, $F'_{d'-2}$  produces a syzygy of $f{\bf g},g$ of degree $d'-1$  whose the $(n+1)$th coordinate is nonzero. Thus, by Proposition~\ref{resolution_of_deJonq}, $d'\geq  d.$ Hence, $d=d'.$ Thus, the equality $\mathfrak{g}'=\mathfrak{g}/\mathfrak{h}$ implies $\mathfrak{h}=1.$ Therefore, the Claim follows.
\qed

\smallskip

Let again $f,g\in R=k[x_1,\ldots,x_{n},x_{n+1}]$ be relatively prime forms  of degree $d-1$ and $d\geq 2$, respectively and set $I:=(x_1f,\ldots,x_{n}f,g)$.

By Corollary~\ref{saturation_CM}, $(x_1,\ldots,x_{n})$ is an associated prime of  $R/I,$  while $\fm=(x_1,\ldots,x_{n+1})$  is not.
Therefore, $R$ admits a linear form  of the shape   $\ell=x_{n+1}-\lambda,$ regular on $R/I$, where $\lambda\in k[x_1,\ldots,x_{n}].$ 
Then the $k$-algebra surjection
$$R=k[x_1,\ldots,x_{n+1}]\surjects k[x_1,\ldots,x_{n}],\quad x_i\mapsto x_i\, (1\leq i\leq n), \,\, x_{n+1}\mapsto \lambda$$
yields an identification $R/(\ell)\simeq S:=k[x_1,\ldots,x_{n}].$

Now, set $J=(I,\ell)/(\ell)\subset R/(\ell)=S.$ Since $\ell$ is regular both on $R$ and on $R/I,$ the ring $S/J$ is a specialization of $R/I.$ In particular, $J$ is an ideal of height 2 minimally generated by $n+1$ form of degree $d$ whose syzygy matrix is the corresponding specialization of $\psi.$

Denote by $\mathcal{J}\subset S[y_1,\ldots,y_{n+1}]$ the presentation ideal of the Rees algebra $\mathcal{R}_S(J).$ One has $(\mathcal{I},\ell)/(\ell)$ upon the identification $R[y_1,\ldots,y_{n+1}]/(\ell)=S[y_1,\ldots,y_{n+1}].$ For a form $P\in R[y_1,\ldots,y_{n+1}]$ denote by $\overline{P}$ its image in $S[y_1,\ldots,y_{n+1}].$ 

Drawing upon the images of the downgraded sequence $F_0,F_1,\ldots,F_{d-2}$, one knows by  \cite[Theorem 3.1]{BenDandr} that  $\mathcal{J}=(\overline{p}_{i,j},\, 1\leq i<j\leq n,\,\overline{F}_0,\overline{F}_1,\ldots,\overline{F}_{d-2}, \mathfrak{H})$ where $\mathfrak{H}$ is a form of bidegree $(0,d)$ in $S[y_1,\ldots,y_{n+1}]$.

\begin{Theorem}\label{main_Rees_structure}
	Let  $I\subset R=k[x_1,\ldots,x_n,x_{n+1}]$ denote the base ideal of a generalized de Jonqui\`eres map $\mathfrak{J}$  with identity support and let
	$\mathcal{I}\subset R[y_1,\ldots,y_{n+1}]$ be the presentation ideal of the Rees algebra $\mathcal{R}_R(I)$ of $I$ obtained by mapping
	$y_i\mapsto x_if$, for $1\leq i\leq n$, and $y_{n+1}\mapsto g$. Then
	\begin{equation}\label{gens_of_ideal_of_Rees_deJonq}
		\{p_{i,j},\, 1\leq i<j\leq n,\,\,F_0,F_1,\ldots,F_{d-2}\}
	\end{equation}
	is a minimal set of generators of $\mathcal{I}$.
\end{Theorem}
\demo
Let $P\in \mathcal{I}_{(i,j)}$ be a form of bidegree $(i,j).$ Consider two cases:

($j< d$) Induct on $i.$ 

There is nothing to prove if $i=0$ since  $\mathfrak{J}$ is in particular a dominant map. Let $i\geq 1.$ Then $\overline{P}\in \mathcal{J}_{(i,j)}$ and, since $j\leq d,$ it cannot involve effectively the form $\mathfrak{H}$, hence we can write
$$\overline{P}=\sum_{1\leq u<v\leq n}\overline{A}_{u,v}\overline{p}_{u,v}+\overline{B}\overline{F}_0+\sum_{u=1}^{d-2}\overline{C}_u\overline{F}_{u}$$
for certain forms $\overline{A}_{u,v},\overline{B}, \overline{C}_u $ of appropriate bidegrees in $S[y_1,\ldots,y_{n+1}].$ Lifting back to $R[y_1,\ldots,y_{n+1}]$ yields
$$P-\sum_{1\leq u<v\leq n} A_{u,v}p_{u,v}-BF_0- \sum_{u=1}^{d-2}C_uF_{u}=\ell Q$$
for a form $Q\in R[y_1,\ldots,y_{n+1}]$ of bidegree $(i-1,j).$ Clearly, then $\ell Q\in \mathcal{I}_{(i,j)},$ hence $Q\in\mathcal{I}_{i-1,j}$ since $\mathcal{I}$ is prime. Thus, by induction $Q\in (p_{i,j},\, 1\leq i<j\leq n,\,\,F,F_1,\ldots,F_{d-2}).$ Therefore $P\in (p_{i,j},\, 1\leq i<j\leq n,\,\,F_0,F_1,\ldots,F_{d-2}).$

\medskip

($j\geq d$) Let $\mathfrak{J}^{-1}=(\mathfrak{g}_1:\ldots:\mathfrak{g}_{n+1}):\pp^n\dasharrow\pp^n$ be the inverse to $\mathfrak{J}.$ By Proposition~\ref{inverse_of_deJonq}, $\mathfrak{g}_i$ is a form in $k[y_1,\ldots,y_{n+1}]$ of degree $d.$

{\sc Claim.} $\ell(\mathfrak{g}_1,\ldots,\mathfrak{g}_{n+1})$ and $\mathfrak{H}$ coincide up to a scalar multiplier. 

To see this, let $\delta\in R$ denote the source inversion factor of $\mathfrak{J}.$ That is, $\delta$ is a homogeneous polynomial of $R$ of degree $d^2-1$  such that $\mathfrak{g}_i(x_1f,\ldots,x_{n}f,g)=x_i\delta,\quad(1\leq i\leq n+1).$
Hence, $\ell(\mathfrak{g}_1,\ldots,\mathfrak{g}_{n+1})(x_1f,\ldots,x_{n}f,g)=\ell(\delta x_1,\ldots,\delta x_{n+1})=\delta\ell(x_1,\ldots,x_{n+1}).$
Reducing modulo $\ell$ we have
$$\ell(\mathfrak{g}_1,\ldots,\mathfrak{g}_{n+1})(\overline{x_1f},\ldots,\overline{x_{n}f},\overline{g})=0.$$ This means that $\ell(\mathfrak{g}_1,\ldots,\mathfrak{g}_{n+1})$ is a multiple of the implicit equation $\mathfrak{H}.$ But $\deg \mathfrak{H}=d=\deg\ell(\mathfrak{g}_1,\ldots,\mathfrak{g}_{n+1}).$ Hence, $\ell(\mathfrak{g}_1,\ldots,\mathfrak{g}_{n+1})$ and $\mathfrak{H}$ are the same up to a scalar multiplier as required. 

Again induct on $i.$ As before, we can assume that $i\geq 0.$ Since $j\geq d$ then 
$$\overline{P}=\sum_{1\leq u<v\leq n}\overline{A}_{u,v}\overline{p}_{u,v}+\overline{B}\overline{F}_0+\sum_{u=1}^{d-2}\overline{C}_u\overline{F}_{u}+\overline{D}\ell(\mathfrak{g}_1,\ldots,\mathfrak{g}_{n+1}),$$
for appropriate coefficients. Hence,
$$P-\sum_{1\leq u<v\leq n}A_{u,v}p_{u,v}-BF_0-\sum_{u=1}^{d-2}C_uF_{u}=D\ell(\mathfrak{g}_1,\ldots,\mathfrak{g}_{n+1})+\ell Q.$$
for some $Q\in R[y_1,\ldots,y_{n+1}]$ of bidegree $(i-1,j).$ Then,
\begin{equation}\label{hip_ind}
	D\ell(\mathfrak{g}_1,\ldots,\mathfrak{g}_{n+1})+\ell(x_1,\ldots,x_{n+1})Q\in\mathcal{I}_{i,j}.
\end{equation}

Now, letting $I':=(\mathfrak{g}_1,\ldots,\mathfrak{g}_{n+1}),$ we know that  there is a $k$-algebra isomorphism $\mathcal{R}_R(I)\simeq \mathcal{R}_R(I')$ induced by the identity map of the bigraded polynomial ring $R[y_1,\ldots,y_{n+1}]=k[x_1,\ldots,x_{n+1},y_1,\ldots,y_{n+1}]$ (\cite[Proposition 2.1]{Si_Cremona}). There follows the equality $\mathcal{I}=\mathcal{I}'$ of the respective presentation ideals. Thus, $$D\ell(\mathfrak{g}_1\ldots,\mathfrak{g}_{n+1})+\ell(x_1,\ldots,x_{n+1})Q\in\mathcal{I}_{(i,j)}\subset\mathcal{I}'$$ as well. This means that, by evaluating $x_i\mapsto \mathfrak{g}_i$ ($1\leq i\leq n+1$), \eqref{hip_ind} gives
$$D(\mathfrak{g}_1,\ldots,\mathfrak{g}_{n+1},y_1,\ldots,y_{n+1})\ell(\mathfrak{g}_1,\ldots,\mathfrak{g}_{n+1})+\ell(\mathfrak{g}_1,\ldots,\mathfrak{g}_{n+1})Q(\mathfrak{g}_1,\ldots,\mathfrak{g}_{n+1},y_1,\ldots,y_{n+1})=0.$$
Since $\ell\notin \mathcal{I},$ $\ell(\mathfrak{g}_1,\ldots,\mathfrak{g}_{n+1})\neq 0;$ canceling yields

$$D(\mathfrak{g}_1,\ldots,\mathfrak{g}_{n+1},y_1,\ldots,y_{n+1})+Q(\mathfrak{g}_1,\ldots,\mathfrak{g}_{n+1},y_1,\ldots,y_{n+1})=0,$$
which means by the same token that
\begin{equation*}
	\vspace{-7pt}
D+Q=D(x_1,\ldots,x_{n+1},y_1,\ldots,y_{n+1})+Q(x_1,\ldots,x_{n+1},y_1,\ldots,y_{n+1})\in \mathcal{I}.
\end{equation*}
But, $D$ and $Q$ have respective bidegrees $(i,j-d)\neq (i-1,j).$ Hence, $D, Q\in\mathcal{I}.$ By the inductive assumption $Q\in (p_{i,j},\, 1\leq i<j\leq n,\,F_0,F_1,\ldots,F_{d-2}).$ As to $D,$ if $j-d<d$ then by the previous argument $D\in (p_{i,j},\, 1\leq i<j\leq n,\,F_0,F_1,\ldots,F_{d-2}).$ Otherwise, we replace $P$ by $D$ and repeat the above argument. Since the second component of the bidegree goes down at each such step, we are through.

The assertion about minimality follows from  \cite{BenDandr}, where it is proved that $\mathcal{J}$ is minimally generated by the set $\{\overline{p}_{i,j},\, 1\leq i<j\leq n,\,\overline{F}_0,\overline{F}_1,\ldots,\overline{F}_{d-2}, \mathfrak{H}\}$.
\qed

\subsubsection{The depth of the Rees algebra}
We now look more closely at the structure of the Rees algebra $\mathcal{R}_R(I)$ of the base ideal of a generalized de Jonqui\`eres map.
Consider the  matrix $(\xx|\yy)$ as in (\ref{symm_algebra_equations}).
The next result seems to be known, but we have not been able to find a precise reference. 

\begin{Lemma}\label{powers_are_primary}
Let $\fn$ denote the residue ideal in $k[\xx,\yy]/I_2(\xx|\yy)$ of the ideal $(\xx)\subset k[\xx,\yy]$. Then the associated graded ring of $\fn$ is isomorphic to $k[\xx,\yy]/I_2(\xx|\yy)$.
In particular, all powers of $\fn$ are $\fn$-primary.
\end{Lemma}
\demo It suffices to show the main assertion because, as is well-known, the associated graded ring being a domain implies that the symbolic powers are ordinary powers.

Set $A:= k[\xx,\yy]/I_2(\xx|\yy)$. The assertion follows from the

{\sc Claim.} Let $\mathcal{R}_A(\fn)\simeq A[\zz]/\mathcal{L}$ be a presentation of the Rees algebra of $\fn$ over $A$, where $\zz=\{z_1,\ldots,z_n\}$ are indeterminates over $A$. Then $\mathcal{L}$ is the residue of the ideal $I_2(G)\subset k[\xx,\yy,\zz]$, where
$$G=\left(\begin{array}{ccccc}x_1&\cdots&x_n\\y_1&\cdots&y_n\\ z_1&\cdots&z_n\end{array}\right).$$

To prove the claim, as is easily seen, $I_2(G)$ is contained in the inverse image of $\mathcal{L}$ in $k[\xx,\yy,\zz]$.
But the first is a prime ideal of height $(3-2+1)(n-2+1)=2(n-1)$ in $k[\xx,\yy,\zz]$.
Therefore, its height in $A[\zz]$ is $2(n-1)-(n-1)=n-1$, because the height of $I_2(\xx|\yy)$ in $k[\xx,\yy]$ is $n-1$.
But we know that the height of $\mathcal{L}$ is one less the number of generators of $\fn$, hence is also $n-1$.
This triggers the required equality.
\qed

\smallskip

Now, $I_2(\xx|\yy)$ fits in the following increasing inclusion sequence of ideals in the presentation ideal $\mathcal{I}$:
$$\mathcal{P}_0:=(I_2(\xx|\yy)) \subset \mathcal{P}_1:=(\mathcal{P}_0,F_0)\subset \mathcal{P}_2:=(\mathcal{P}_1,F_1)\subset\cdots\subset \mathcal{P}_{d-1}:=(\mathcal{P}_{d-2},F_{d-2}).$$
Here, the round brackets indicate an ideal of $R[y_1,\ldots,y_{n+1}]$, where $R=k[x_1,\ldots,x_{n+1}]$.
Note that, by Theorem~\ref{gens_of_ideal_of_Rees_deJonq},  $\mathcal{P}_{d-1}=\mathcal{I}$. 

The following lemma clarifies the role of these ideals.

\begin{Lemma}\label{Fi_is_not}
	{\rm (i) }	Let $\mathfrak{f}\in \mathcal{I}$ be a form of bidegree $(r,s).$ 
	If $\mathfrak{f}\in (\mathcal{P}_0, (x_1,\ldots,x_n)^r)\subset R[y_1,\ldots,y_{n+1}]$  then $\mathfrak{f}\in \mathcal{P}_0.$ In particular, $F_i\notin (\mathcal{P}_0, (x_1,\ldots,x_n)^{d-i}) R[y_1,\ldots,y_{n+1}]$, for every $1\leq i\leq d-2.$
	
	{\rm (ii) } $\mathcal{P}_0:F_0=\mathcal{P}_0$ and $\mathcal{P}_{i}:F_i=(x_1,\ldots,x_n)$ $(1\leq i\leq d-2).$
\end{Lemma}
\demo (i) Writing $\mathfrak{f}=\mathfrak{f}_1+\mathfrak{f}_2$ where $\mathfrak{f}_1\in\mathcal{P}_0$ and $\mathfrak{f}_2\in (x_1,\ldots,x_n)^{r}$ we have that $\mathfrak{f}_2$ is a form of bidegree $(r,s)$ in $\mathcal{I}$ that depends only on the variables $x_1,\ldots,x_n,y_1,\ldots,y_{n+1}.$ Hence, by Proposition~\ref{gens_of_ideal_of_Rees_deJonq} (b), $\mathfrak{f}_2$ is a Rees relation of the ideal $(y_1,\ldots,y_n)$, hence $\mathfrak{f}_2\in\mathcal{P}_0.$

If $F_i\in (\mathcal{P}_0, (x_1,\ldots,x_n)^{d-i}) R[y_1,\ldots,y_{n+1}]$ then, by the first part, we have $F_i\in \mathcal{P}_0$, which contradicts the fact that, by Theorem~\ref{main_Rees_structure},  $\{p_{i,j},\, 1\leq i<j\leq n,\,\,F_0,F_1,\ldots,F_{d-2}\}$ is a  minimal set of generators of $\mathcal{I}.$

\smallskip

(ii) The first equality is clear since $\mathcal{P}_0$ is a prime ideal and $F_0\notin \mathcal{P}_0$, by Remark~\ref{destaques} (a).

Now, take $1\leq i\leq d-2.$ By Remark~\ref{destaques}(c),  $(x_1,\ldots,x_n)\subset \mathcal{P}_{i}:F_i.$ For the other inclusion, consider  $a\in \mathcal{P}_{i}:F_i.$ Then,
\begin{equation}
	aF_i\in \mathcal{P}_{i}
	\subset (\mathcal{P}_0,\,(x_1,\ldots,x_n)^{d-i})R[y_1,\ldots,y_{n+1}],
\end{equation}
where the inclusion  is a consequence of Remark~\ref{destaques}(b). Thus,  $\overline{a}\overline{F_i}\in (\overline{x}_1,\ldots,\overline{x}_n)^{d-i},$
where ``$-$" denotes the residual class of a element module $\mathcal{P}_0=I_2(\xx|\yy)$. By Lemma~\ref{powers_are_primary}, $(\overline{x}_1,\ldots,\overline{x}_n)^{d-i}$ is $(\overline{x}_1,\ldots,\overline{x}_n)$-primary and, by part (i),  $\overline{F_i}\notin(\overline{x}_1,\ldots,\overline{x}_n)^{d-i}$. Therefore,  $\overline{a}\in (\overline{x}_1,\ldots,\overline{x}_n).$ That is, $a\in (x_1,\ldots,x_n)$ since $\mathcal{P}_0\subset (x_1,\ldots,x_n)$.
\qed

\begin{Proposition} Let $R=k[x_1,\ldots,x_n,x_{n+1}]$, with  $n\geq 2$, and let $I\subset R$ denote the base ideal of  a generalized de Jonqui\`eres transformation of degree $d\geq 2$ with identity support. 
	Then the Rees algebra $\mathcal{R}(I)$ is almost Cohen--Macaulay.
\end{Proposition}
\demo 
Set $A=R[y_1,\ldots,y_{n+1}].$ 
We depart from two well-known free resolutions.
The first is the $2$-linear resolution of $A/\mathcal{P}_0$:
\vspace{-7pt}
\begin{equation}\label{resolution_of_2x2}
	0\to G_{n-1,0}\to \cdots\to G_{2,0}\to G_{1,0}\to A\to A/\mathcal{P}_0\to 0,
	\vspace{-5pt}
\end{equation}
where $G_{i,0}=A(-(i+1))^{i\,{n\choose i+1}}$, for $i\geq 1$ and $G_{0,0}=A$.  
The second is the shifted Koszul complex resolving the complete intersection $(A/(x_1,\ldots,x_n))(-d)$:
{\small 
	\begin{equation}\label{resolutionxx}
		0\to K_n(-(d+n)) \to \cdots \to  K_1(-(d+1))\to A(-d)\to (A/(x_1,\ldots,x_n))(-d)\to 0.
		\vspace{-5pt}
	\end{equation}
}
To start, consider the exact sequence $A/\mathcal{P}_0\stackrel{\cdot F_0}{\lar}  A/\mathcal{P}_0 \lar A/(\mathcal{P}_0,F_0)=A/\mathcal{P}_1\rar 0,$
of multiplication by $F_0$ on $A/\mathcal{P}_0$,
using the definition of $\mathcal{P}_1$.
This map extends to a map of (\ref{resolution_of_2x2}) to itself and the resulting mapping cone (\cite[Section 6.2.3.2]{SimisBook}) has the form
\begin{equation}\label{resolutionp1}
	0\to G_{n,1}\to G_{n-1,1}\to\cdots\to G_{1,1}\to  A\to A/\mathcal{P}_1\to0,
\end{equation}
where $G_{i,1}=G_{i,0} \oplus G_{i-1,0}$, for $0\leq i \leq  n.$
This a graded free resolution of $A/\mathcal{P}_1$ of length $n$ with $G_{n,1}=G_{n-1,0}=A(-n)^{n-1}$.

Next, consider the homogeneous map 
$$(A/(x_1,\ldots,x_n))(-d)= (A/\mathcal{P}_1:F_1)(-d) \to A/\mathcal{P}_1$$
induced by multiplication by $F_1$, with cokernel $A/(\mathcal{P}_1, F_1)=A/\mathcal{P}_2$.
It extends to a map from (\ref{resolutionxx}) to (\ref{resolutionp1}), with resulting mapping cone given by
\begin{equation}\label{resolutionp2}
	0\to G_{n+1,2}\to G_{n,2}\to\cdots\to G_{1,2}\to A\to A/\mathcal{P}_2\to0,
\end{equation}
where $G_{i,2}=G_{i,1}\oplus K_{i-1}(-(d+i-1))$. This resolution has length $n+1$. 

Repeat the last argument with multiplication by $F_2$ on $A/\mathcal{P}_2$ and the extended map from (\ref{resolutionxx}) to (\ref{resolutionp2}) above to get a graded free resolution of $A/\mathcal{P}_3$ of length $n+1$.  

Proceeding this way, we eventually reach a graded free resolution of $A/\mathcal{P}_{d-1}=A/\mathcal{I}$ of length $n+1$.
Therefore, ${\rm proj.dim\,} A/\mathcal{I}\leq n+1,$ that is, $\depth  A/\mathcal{I}\geq n+1$, showing that $\mathcal{R}_R(I)$ is almost Cohen--Macaulay. 
\qed

\begin{Conjecture}
With the same notation as in the previous theorem,  $\mathcal{R}_R(I)$ is Cohen--Macaulay if and only if $d\leq n+1$.
\end{Conjecture}
Note that if $\mathcal{R}(I)$ is Cohen--Macaulay then the relation type of $I$ is at most the analytic spread of $I$ (see, e.g., \cite[Corollary 6.3]{AHT}). Therefore, $d-1\leq n+1$, i.e., $d\leq n+2$, it remaining to exclude the case where $d=n+2$. 
For the converse, one follows the same line of proof above, thereby arguing that, if $d\leq n+1$, there is a canceling at the tail of (\ref{resolutionp2}) and the subsequent iterated mapping cones, resulting a minimal resolution of length $n$ at the end.



\noindent {\bf Addresses:}

\medskip

\noindent {\sc Zaqueu Ramos}\\
Departamento de Matem\'atica, CCET\\ 
Universidade Federal de Sergipe\\
49100-000 S\~ao Cristov\~ao, Sergipe, Brazil\\
{\em e-mail}: zaqueu@mat.ufs.br\\

\medskip

\noindent {\sc Aron Simis}\\
Departamento de Matem\'atica, CCEN\\ 
Universidade Federal de Pernambuco\\ 
50740-560 Recife, PE, Brazil\\
{\em e-mail}:  aron@dmat.ufpe.br

\end{document}